\documentclass[12pt,a4paper]{amsart}

\usepackage{amsmath, amsfonts, xifthen, latexsym, amssymb, amsthm, amscd}

\usepackage[utf8]{inputenc}
\usepackage[shortlabels]{enumitem}
\usepackage{graphicx}
\usepackage{url}

\usepackage{hyperref}
\hypersetup{colorlinks=true,citecolor=blue,filecolor=blue,linkcolor=blue,urlcolor=blue}
\usepackage[margin=1.4in]{geometry}

\newtheorem{theorem}{Theorem}

\newtheorem{prop}{Proposition}

\theoremstyle{definition}
\newtheorem{definition}{Definition}

\newcommand{\eps}{{\varepsilon}}

\renewcommand{\phi}{{\varphi}}

 \newcommand\spec{\operatorname{Spec}}

\newcommand\diam{\operatorname{diam}}

\newcommand{\cB}{\mathcal{B}}
\newcommand{\cCB}{\mathcal{CB}}

\newcommand\N{\mathbb N}
\newcommand\cG{\mathcal{G}}

\newcommand{\actson}{\curvearrowright}

\newcommand{\cC}{\mathcal{C}}
\newcommand{\cant}{\{0,1\}^\N}
\newcommand{\cal}[1]{{\mathcal #1}}

\usepackage{etoolbox}
\newtoggle{no_cases}\toggletrue{no_cases}
\newcommand{\case}[2][]{\iftoggle{no_cases}{\left\{\begin{array}{ll}#2 & #1}{\\#2 & #1}\togglefalse{no_cases}}
\newcommand{\esac}{\end{array}\right.\toggletrue{no_cases}}

\newcommand{\Prob}{\operatorname{Prob}}
\newcommand{\stab}{\operatorname{Stab}}\newcommand{\Sub}{\operatorname{Sub}}

\newcommand{\gd}{FG_d}
\newcommand{\vi}{\vskip 0.1in \noindent}
\begin{document}
\title[Cantor combinatorics]{Cantor combinatorics and almost finiteness}
\author{Gábor Elek}
\address{Department of Mathematics And Statistics, Fylde College, Lancaster University, Lancaster, LA1 4YF, United Kingdom}

\email{g.elek@lancaster.ac.uk}  

\begin{abstract} In this survey we give a concise introduction to a continuous version of 
Borel
combinatorics based on \cite{Elekcantor} and \cite{Elekurs}.
Our approach has a certain algorithm-theoretic nature and we are giving special emphasis to the notion of almost finiteness introduced by Matui \cite{Matui}
as a continuous
analogue of Borel hyperfiniteness. 

\end{abstract}\maketitle
\vi
\vi
{\bf 1. The Properness Assumption.} \vi
In the last couple of decades Borel combinatorics has become
an important branch of mathematics (see e.g. \cite{kema} and \cite{kemi} for surveys). Let us briefly recall the most
important definitions. Let $X$ be the standard Borel space. A {\bf Borel graph}
is a subset $E\subset X\times X$ such that if $(x,y)\in E$ then $(y,x)\in E$. 
Also, if $x\in X$, then $(x,x)\in E$ implies
that $(x,y)\notin E$ for any $x\neq y\in E$ (so, we do not allow loops only
singletons).
 We will assume the {\bf bounded degree}
condition as well. There exists some integer $d>0$ such that if $x\in X$, then
$$\deg(x)=|\{y\in X\,\mid\, (x,y)\in E\}|\leq d\,.$$
\noindent
It is easy to see, that the component of such Borel graphs $E$ are always
finite or countably infinite graphs. Let $\Gamma$ be a finitely generated group
with symmetric generating system $\Sigma$ and let $\alpha:\Gamma\actson X$ be
an action of the group $\Gamma$ on $X$ by Borel isomorphisms.
If $x\neq y$, let $(x,y)\in E_\alpha$ if $\alpha(\sigma)(x)=y$ for some $\sigma\in\Sigma$.
Also, let $(x,x)\in E_\alpha$ if $\alpha(\gamma)(x)=x$ for
all $\gamma\in\Gamma$. Then, $E_\alpha\subset X\times X$ is a Borel graph and it is
called the graph of the action $\alpha$ with respect to $\Sigma$.

\noindent
Let us mention a basic result on
Borel graphs. 
\begin{prop} [Kechris-Solecki-Todorcevic,\cite{KST}] \label{Brooks} Every Borel graph
of vertex degree bound $d$ can be properly colored by $d+1$ colors.
That is, there exists a partition of $X$ into $d+1$ Borel sets
$A_1, A_2,\dots, A_{d+1}$ such that if $x\in A_i$ and $(x,y)\in E$, then
$y\notin A_i$.
\end{prop}
\noindent
As a corollary one can immediately see, that if $E$ is a Borel graph of vertex degree bound $d$, then
there exists a Borel action $\alpha: \Gamma \actson X$ such that $E_{\alpha}=E$. In fact, $\Gamma$ can be chosen as the
free product of $\binom{d}{2}$ copies of the cyclic group of two elements.

\noindent
One of the most important notions in the theory of Borel graphs is {\bf hyperfiniteness}.
A Borel graph is called hyperfinite if there exist Borel subgraphs
$E_1\subset E_2\subset\dots$ such that $\cup_{i=1}^\infty E_i=E$ and all
the components of the graphs $E_i$ are finite.
It is conjectured that Borel graphs of free amenable actions are
always hyperfinite. 
\begin{prop}[Jackson-Kechris-Louveau \cite{JKL}] \label{poly} Let $\Gamma$ be
a finitely generated group of polynomial growth. Then for any Borel action
of $\Gamma$ on a standard Borel set, the associated Borel graph is hyperfinite.
\end{prop}
\noindent
The goal of this note is to introduce the Cantor version of 
Borel graph combinatorics and we will use the two propositions above to test our
theory. 
\begin{definition} A {\bf Cantor graph} is a 
Borel graph $E\subset \cC\times\cC$, 
such that $E$ is a closed subset, where $\cC=\cant$ is the Cantor set.
\end{definition}
\noindent
Let us observe that the Cantor analogue of Proposition \ref{Brooks} does not
hold. Consider a 
continuous action $\alpha:\Gamma\actson \cC$  of a finitely generated group $\Gamma$ with generating system 
$\Sigma$,
such that for a certain $x\in \cC$ and for some generator
$\sigma\in\Sigma$, $\alpha(\sigma)(x)=x$ however, there exists a sequence $\{x_n\}_{n=1}^\infty\to x$ such that $\alpha(\sigma)(x_n)=y_n\neq x_n$. 
The simplest example for such systems is the standard Bernoulli action $\alpha:\Gamma\actson \{0,1\}^\Gamma$, where $x\in \{0,1\}^\Gamma$ is constant
as a $\{0,1\}$-valued function on $\Gamma$. 
Then, in the Borel graph $E_\alpha$ associated to the action $\alpha$, $x_n$ is adjacent to $y_n$ for all $n\geq 1$. Suppose that
we have a clopen partition $\cC=\cup^q_{i=1} A_i$ and $x\in A_i$. Then for large enough $n$, both $x_n$ and $y_n$ are in $A_i$, hence
the partition does not define a clopen proper coloring of the Borel graph $E_\alpha$.

\noindent
In order to avoid the obstacle above, we introduce the notion of {\bf proper} Cantor graphs. Note that if $E$ is a Cantor graph and
$x,y\in \cC$, then $d_E(x,y)=r$ if $x$ and $y$ are in the same component of $E$ and their shortest path distance is $r$. Otherwise, let
$d_E(x,y)=\infty$. Also, if $x\in \cC=\cant$ and $t$ is a positive integer,  then $(x)_t\in \{0,1\}^{[t]}$ is the projection of $x$ onto
its first $t$ coordinates.  
\begin{definition}\label{propercantor}
The Cantor graph $E$ is called proper if for any $r>1$, there exists $s_r>1$ such that if $d_E(x,y)\leq r$ then $(x)_{s_r}\neq (y)_{s_r}$.
\end{definition}
\noindent
We will prove that following Cantor analogue of Proposition \ref{Brooks}.
\begin{prop} \label{topKST}
Any proper Cantor graph of vertex degree bound $d$ can be properly colored by $d+1$ colors in a continuous way.
That is, there exists a partition of $\cC$ into $d+1$ clopen sets
$A_1, A_2,\dots, A_{d+1}$ such that if $x\in A_i$ and $(x,y)\in E$, then
$y\notin A_i$.
\end{prop}
\noindent
Consequently, any proper Cantor graph can be obtained as a graph of a continuous Cantor action of a finitely generated group. 
\noindent
Recall \cite{Elekurs} that a continuous action of a finitely group $\Gamma$, $\alpha:\Gamma\actson \cC$ is
called  proper if all the stabilizer groups
are stable. That is,  if $\alpha(\gamma)(x)=x$ for some $\gamma\in\Gamma$  then there exists
some open set $x\in U$, such that $\alpha(\gamma)(y)=y$ provided that $y\in U$. Clearly, the associated Cantor graph of a proper action is proper. Of course, the simplest examples for proper actions are the free actions.

\noindent
The question is, whether there are many interesting {\bf non-free} proper actions (and hence proper Cantor graphs) as well. The notion of  uniformly recurrent
subroups, introduced by Glasner and Weiss \cite{GW}, provides a hopefully satisfactory ``yes'' answer for this question. 
If an action $\alpha$ is proper, then the stabilizer map
$\stab_\alpha:\cC\to \Sub(\Gamma)$ is a continuous map from the Cantor set to the compact space of all subgroups of $\Gamma$. This map is equivariant
with respect to $\alpha$ and to the conjugacy action on $\Sub(\Gamma)$.
The action is free if the image of $\stab_\alpha$ contains only the trivial subgroup. In general, the image of $\stab_\alpha$ for a proper
action is a closed, invariant subspace of $\Sub(\Gamma)$. A {\bf uniformly recurrent subgroup} $Z\subset \Sub(\Gamma)$ is a closed, minimal invariant subspace
(that is, all orbits are dense). Answering a query
of Glasner and Weiss \cite{GW}, it is proved in \cite{Elekurs} that if $Z$ is a uniformly recurrent
subgroup in $\Sub(\Gamma)$ then there is a proper action $\alpha$ so that the image of $\stab_\alpha$ is just $Z$.
This means that there is huge zoo of proper minimal Cantor graphs with interesting properties. For an example, if $\beta>1$ is a real number then there exists a proper minimal
Cantor graph $E$ such that all of the growth exponent of its orbit is $\beta$. Note that a Cantor graph
is minimal if there exists no proper compact subset $T\subset \cC$ such that $T$ is the union of components.
\vi
{\bf 2. Almost Finiteness.} \vi
Now we define the Cantor version of almost finiteness following Matui \cite{Matui} and Suzuki \cite{Suzuki}.
Let $\alpha:\Gamma\actson \cC$ be a continuous Cantor action as above, where $\Gamma$ is a finitely generated group with a symmetric generating set $\Sigma$.
A {\bf clopen tower} $(S,F)$ of $\alpha$ consists of a clopen subset $S\subset \cC$ and a finite subset $F\subset\Gamma$ such that
if $s\neq t\in S$, then $\alpha(F)(s)\cap\alpha(F)(t)=\emptyset.$ We call $\alpha(F)(s)$ a {\bf tile} of the tower. 
Recall that if $G$ is a graph of bounded vertex degrees and $H\subset V(G)$ is a finite subset then the {\bf isoperimetric constant} of $H$
is defined as
$$i_G(H)=\frac{|\partial(H)|}{|V(H)|}\,,$$
where $x\in\partial(H)$ if $x\in H$ and there exists
$y\in V(G)\backslash H$ such that
$y$ is adjacent to $x$.
The isoperimetric constant of a tower $(S,F)$ is defined as
$$i_{\alpha}((S,F))=\sup_{s\in S} i_{G_s}(\alpha(F)(s))\,,$$
where $G_s$ is the component of $s$ in the graph $E_\alpha$.
A {\bf full clopen $\eps$-castle} of $\alpha$, $\cup^m_{j=1} (S_j,F_j)$ is a union of towers such that
\begin{itemize}
\item $\cup^m_{j=1} S_j$ is a partition of $\cC$.
\item For any $1\leq j \leq m$, $i_\alpha(S_j,F_j)\leq \eps\,.$
\end{itemize}
\begin{definition}
The action $\alpha:\Gamma\actson \cC$ is {\bf almost finite } if for any $\eps>0$, $\alpha$ admits
a full clopen $\eps$-castle. A proper Cantor graph $E$ is almost finite if there is an almost finite action $\alpha$ such that $E_\alpha=E$.
\end{definition}
\noindent
In order to state our main result on almost finiteness we need to recall the notion of $D$-doubling.
Let $G$ be a graph of bounded vertex degrees and $D>0$ be an integer. We say that $G$ is $D$-doubling if for any $x\in V(G)$ and $s\geq 1$, we have that
$$|B_{2s}(G,x)|\leq D |B_s(G,x)|\,,$$
\noindent
where $B_s(G,x)$ is the ball of radius $s$ around the vertex $x$ in the graph $G$.
Note that the definition above can be applied to finite graphs, countable graphs and Cantor graphs as well. 
Although it is not true that all the graphs of polynomial growth are doubling, graphs
of strict polynomial growth are always doubling. Recall that a graph $G$ is of
strict polynomial growth if there exists $\alpha\geq 1$, $C_1, C_2>0$ such that
$$C_1 r^\alpha\leq B_r(G,x) \leq C_2 r^\alpha$$
holds for all $x\in V(G)$ and $r\geq 1$. Therefore, the following theorem can be viewed as a Cantor analogue of Proposition \ref{poly}.
\begin{theorem} \label{fotetel} Let $E\subset \cC\times \cC$ be a $D$-doubling proper Cantor graph. Then $E$ is almost finite. 
\end{theorem}
\noindent
In the proof of Theorem \ref{fotetel} we use a ``non-free'' version of the idea of  Downarowicz, Huczek and Zhang
\cite{DHZ}, that does not seem to work in
the original case of amenable groups, but works nicely when the graphs are doubling. Note that Theorem \ref{fotetel} implies that all the free actions
of a finitely generated nilpotent-by-finite group is almost finite. It is an open question whether continuous actions of amenable groups are always almost finite.
Nevertheless, modifying the construction of Section 10 in \cite{Elekurs} we will prove the following result.
\begin{theorem}
There exist proper, minimal Cantor actions that are almost finite but not topologically amenable 
(see \cite{Oza} for a survey on topological amenability).
\end{theorem}
\noindent
One can define the almost finiteness of a countable graph $G$ the following way.
\begin{definition} A countable graph $G$ of bounded vertex degrees is called almost finite if for any $\eps>0$ there exists $K_\eps>0$ and a partition of
$V(G)$ into subsets $L_1, L_2,\dots$ such that
\begin{itemize}
\item For any $i\geq 1$, $\diam(L_i)\leq K_\eps$.
\item $i_G(L_i)\leq \eps$. 
\end{itemize}
\end{definition}
\noindent
The main result of \cite{DHZ} is that the Cayley graph of a finitely generated amenable group is always almost finite. We can prove that
any doubling countable graph is almost finite. 
\vi {\bf 3.  Fractional Almost Finiteness.}\vi
The notion of fractional almost finiteness is motivated by the notion of fractional hyperfiniteness introduced by L\'aszl\'o Lov\'asz \cite{Lov}. Fractional
almost finiteness will be crucial in the proof of Theorem \ref{specconv} and of Theorem \ref{luck}, but seems to be an interesting notion on its own.
\begin{definition}
A graph $G$ of bounded vertex degrees is called {\bf fractionally almost finite} if for 
any $\eps>0$ and $0<p<1$ there exists
a $K>0,Q>0$ and partitions $\{\cup_{i=1}^\infty L_i^n\}^Q_{n=1}=V(G)$ such that 
\begin{itemize}
\item For any $i\geq 1$ and $1\leq n \leq Q$, $\diam(L_i^n)\leq K$.
\item For any $x\in V(G)$, 
$$|\{n\,\mid\, x\in L^n_i\backslash\partial(L^n_i)\,\mbox{for some $i\geq 1$}\}|\geq pQ\,.$$
\end{itemize}
\end{definition}
\noindent
That is, all the vertices are ``mostly'' inside the tile they are contained in, not on the boundary 
of the tile.
Clearly, the Cayley graph of a non-amenable group cannot be almost finite, so the following result might be surprising.
\begin{theorem}
The Cayley graph of a free group is fractionally almost finite. On the other hand, the Cayley-graphs of non-exact groups
such as the Gromov-Osajda groups are not fractionally almost finite
\end{theorem}
\vi
{\bf 4. Benjamini-Schramm convergence vs. topological convergence.}
\vi
Let us recall one of most important notions of measurable graph limit theory: 
Benjamini-Schramm convergence (see e.g. the monograph of L\'aszl\'o Lov\'asz \cite{Lov2})
Let $d>0$ and $FG_d$ denote the set of all connected, simple finite graphs (up to isomorphism) with vertex degree bound $d>0$.
Let $U^k_d$ be the set of all rooted balls of radius at most $k$ with vertex degree bound $d$. That is, an element of
$U^k_d$ is a graph $H\in \gd$ with a distinguished vertex $x$, such that if $y\in V(H)$, then $d_H(x,y)\leq k$, where
$d_H$ is the shortest graph distance on the graph $H$. For a finite graph $G\in FG_d$ we have a natural probability distribution
$\Prob^k_G$ on $U^k_d$ defined in the following way. If $B\in U^k_d$, then
$$\Prob^k_G(B)=\frac{|P(G,B)|}{|V(G)|}\,,$$
where $P(G,B)$ is the set of vertices $y\in V(G)$ such that the ball of radius $k$ around $y$, $B_k(G,y)$ is rooted-isomorphic to
$B$. 
Let $\{G_n\}^\infty_{n=1}\subset \gd$ be a sequence of finite graphs. We say that  $\{G_n\}^\infty_{n=1}$ is {\bf convergent in the
sense of Benjamini and Schramm} if for any $k\geq 1$ and $B\in U^k_d$, $\lim_{n\to\infty} \Prob^k_{G_n}(B)$ exists.
Now let $E\subset X\times X$ be a Borel graph with an invariant measure $\mu$ (that is, $\mu$ is invariant under any group action
that defines $E$). Then, for each $k\geq 1$ and $B\in U^k_d$ we have a Borel set $V_B\subset X$ such that $B_k(E,x)$ is rooted-isomorphic
to $B$ if and only if $x\in V_B$. We say that $(E,\mu)$ is a Benjamini-Schramm limit of the convergent sequence of finite graphs $\{G_n\}^\infty_{n=1}$
if for any ball $B$
$$\lim_{n\to\infty} \Prob^k_{G_n} (B)=\mu(V_B)\,.$$
\noindent
For every Benjamini-Schramm convergent graph sequence one can find a Benjamini-Schramm limit and it is not known whether 
there exists a convergent graph sequence for any pair $(E,\mu)$.

\noindent
Now we define the topological analogue of the Benjamini-Schramm convergence.
Let $Gr_d$ denote the set of all (not necessarily connected) graphs of vertex degree bound $d$.
For $G\in Gr_d$ denote by $\cB(G)$ the set of rooted balls $B$ for which there is an $x\in V(G)$ and $k\geq 1$
such that $B_k(G,x)$ is rooted isomorphic to $B$.
Let $G,H\in Gr_d$. We say that
$G$ and $H$ is {\bf equivalent}, if $\cB(G)=\cB(H)$. 
We will denote by $\overline{Gr_d}$ the set of equivalence classes of $Gr_d$.
We can define a metric on $\overline{Gr_d}$ as follows.
Let $G,H\in Gr_d$ representing the elements $[G],[H]\in \overline{Gr_d}$. Then
$d_{Gr}([G],[H])= 2^{-n}$ if
\begin{itemize}
\item
For any $1\leq i \leq n$ and ball $B\in U^i_d$, either $B\in \cB(G)$ and $B\in \cB(H)$, 
or $B\notin \cB(G)$ and $B\notin \cB(H)$.
\item There exists $B\in U^{n+1}_d$ such that $B$ is a rooted ball in exactly one of the two graphs.
\end{itemize}
\noindent
Then the space $\overline{Gr_d}$ is compact with respect to the metric above. We call convergence in the metric $d_{Gr}$ topological convergence.
We can also define a Cantor graph limit of topologically convergent graph sequences similarly to the Benjamini-Schramm limit.
Let $G$ be a countable graph, we say that a proper Cantor graph $E\subset \cC\times \cC$ is a realization of $G$ if $\cB(G)$ equals to the set of rooted
balls in $E$. It is not hard to see that for  each $G$ one can find a proper Cantor graph realization $E_G$ and conversely, any proper Cantor graph
$E$ is a realization of a countable graph $G$.  If $\{G_n\}^\infty_{n=1}$ topologically converges to $G$, then we say that $E_G$ is
a Cantor graph limit of $\{G_n\}^\infty_{n=1}$. One can show that if $G$ is the topological limit of a finite graph sequence $\{G_n\}^\infty_{n=1}$, then
it has a realization that admits an invariant measure. Nevertheless, we have the following proposition. 
\begin{prop}
There exist graphs $G$ such that none of the realizations of $G$ admit invariant measure. Also, there exists $G$ such that all of its realizations
admit an invariant measure, but $G$ is not the topological limit of finite graph sequences.
\end{prop}
\noindent
We can define the topological limit of Cantor labeled graph sequences as well, and these objects are crucial in the study of proper Cantor graphs.
A Cantor labeled graph is a countable graph $G\in Gr_d$ equipped with a vertex labeling
$\phi:V(G)\to \{0,1\}^\N$. The set of Cantor labeled graphs is denoted by $CGr_d$.
For $k\geq 1$ we denote by $CU^k_d$ the set of rooted balls $B$ of radius $d$ equipped with a vertex labeling
$\rho:V(B)\to \{0,1\}^{[k]}$, where $[k]=\{1,2,\dots,k\}$. So, if $(G,\phi)\in CGr_d$ and $k\geq 1$, then for any $x\in V(G)$ we assign
an element of $CU^k_d$. Namely, the the $k$-ball centered at $x$ labeled by $\phi_k$, where 
$\phi_k(z)=(\phi(z))_k$ (the projection onto
the first $k$ coordinates).  Now we can proceed exactly the same way as in the unlabeled case. For $G\in CGr_d$
and $B\in CU^k_d$, $B\in \cCB(G)$ if and only if there exists $x\in V(G)$ such that the rooted $\{0,1\}^{[k]}$-labeled ball $B_k(G,x,\phi_k)$ is rooted-labeled
isomorphic to $B$.
We say that $(G,\phi)$ and $(H,\psi)$ are equivalent if $\cCB(G)=\cCB(H)$. The set of
equivalence classes will be denoted by $\overline{CGr_d}$. The metric
on $\overline{CGr_d}$ is defined as follows. Let $(G,\phi),(H,\psi)\in CGr_d$ representing
the classes $[(G,\phi)]$ and $[(H,\psi)]$. Then,
$$d_{CGr}((G,\phi),(H,\psi))=2^{-n}\,,$$
if
\begin{itemize}
\item For any $1\leq i \leq n$ and $B\in CU^i_d$, $B\in \cCB(G)$, either $B\in \cCB(G)$ and $B\in \cCB(H)$, 
or $B\notin \cCB(G)$ and $B\notin \cCB(H)$.
\item There exists $B\in CU^{n+1}_d$ such that
$B$ is rooted-labeled isomorphic to a $\{0,1\}^{[n+1]}$-labeled ball of exactly one of the two graphs.
\end{itemize}
\noindent
Again, the metric space $(\overline{CGr_d}, K)$ is compact. Now let $E\in \cC\times \cC$ be a Cantor graph of vertex degree bound $d$. 
If $x\in\cC$, then the component
of $x$ can be considered as a Cantor labeled graph, where the label of a vertex $y$ is itself. We say that a Cantor graph is the labeled realization
of a Cantor labeled graph $(G,\phi)$ if the  $\{0,1\}^{[k]}$-labeled balls in $E$ are exactly the \\ $\{0,1\}^{[k]}$-labeled balls in $(G,\phi)$. Any proper Cantor graph is the Cantor labeled realization of some $(G,\phi)\in CGr_d$.  
\vi{\bf 5. Almost finiteness vs. Hyperfiniteness.}\vi
First let us recall the notion of hyperfiniteness \cite{Elekhyper}.
A class $\cal{H}\subset FG_d$ of finite graphs is {\bf hyperfinite} if for any
$\eps>0$ there exists $K>0$ such that for any $G\in \cal{H}$ one can delete $\eps |V(G)|$ edges from $E(G)$ in
such a way that in the resulting graph $G'$ all the components have size less or equal than $K$.
The classes of planar graphs or graphs of fixed polynomial growth are hyperfinite.
It is not hard to see that if the class $\cal{H}$ is almost finite (that is the union of the elements
as a countable graph is almost finite), then $\cal{H}$ is hyperfinite as well. On the other hand,
the class of finite trees
is hyperfinite but not almost finite. Note that Matui has the following result:
\cite{Matui}: If $\alpha:\Gamma\actson \cC$ is almost finite, then for any 
$\alpha$-invariant probability measure $\mu$,
the measured equivalence relation defined by $\alpha$ is hyperfinite (see \cite{kemi} for the definition
of measurable hyperfiniteness). On the other hand, it is not known whether almost finiteness of a 
Cantor action implies Borel
hyperfiniteness or not. If $\{G_n\}^\infty_{n=1}$ is a a Benjamini-Schramm convergent sequence of graphs of vertex degree bound $d$, then
all the  Benjamini-Schramm limits of $\{G_n\}^\infty_{n=1}$ are measurably hyperfinite \cite{Schramm}, \cite{Elekhyper} if and only if the family $\{G_n\}^\infty_{n=1}$ is hyperfinite.
It is not hard to see that if $\{G_n\}^\infty_{n=1}$ is a topologically convergent sequence tending to a countable graph $G$, then
$G$ has at least one almost finite proper Cantor graph realization. On the other hand, it is not known
that the existence of such almost finite realization implies the almost finiteness of the sequence.

\vi {\bf 6.  Local algorithms.} \vi
{\bf Local algorithms} (or local distributed algorithms) are algorithms on graphs performed in
parallel by simple machines that have access only small parts of the graphs. In order to illustrate the notion let us consider
the proper vertex coloring problem.
The input is a finite graph of vertex degree bound $d$. At each vertex $x\in G$ there is a  machine
that color the vertex $x$ by one of the elements of the set $Q$, $|Q|=d+1$, by searching the $k$-neighbourhood of $x$,
where $k\geq 1$ does not depend on $x$. The goal is to color adjacent vertices by different colors.
Formally speaking, there is  a function ({\bf the oracle}) $\theta:U^k_d\to Q$ and the coloring $\phi:V(G)\to Q$ is
defined by the rule
$$\phi(x)=\theta(B_k(G,x))\,,$$
\noindent
for each vertex $x\in V(G)$.
Unfortunately, local algorithms can meet obstacles to perform such tasks like the one above due to 
the possible symmetries of the graph $G$.
Say, we have a cycle $C_n$ of length $n$ and we try to color it by three colors.
Since for $k\geq 1$, all the balls $B_k(C_n,k)$ are isomorphic oracles can compute only constant-valued colorings. In order to
break the symmetries, one can use coin tosses to obtain a uniformly random labeling $\lambda:V(G)\to \{0,1\}$. Then the
machine at the vertex $x$ check the labeled ball $B_k(G,x,\lambda)$ before making the decision. One can actually prove that for any $\eps>0$
there exists $k\geq 1$ and $\theta: U^{k,\{0,1\}}_d\to \{1,2,\dots,d+1\}$ such that for any finite graph $G\in G_d$ and
for the resulting vertex coloring $\lambda^\theta:V(G)\to \{1,2,\dots, d+1\}$ the following statement holds:
The probability of the event that the number of vertices $x\in V(G)$ for which the color of $x$ is different from the colors of its 
neighbours is larger than $(1-\eps) |V(G)|$ is greater than $(1-\eps)$ (note that $U^{k,\{0,1\}}_d$ denotes the set of equivalence classes of
rooted $k$-balls of vertices labeled by the set $\{0,1\})$. That is, using randomized local algorithms one can properly color
most of the vertices of a finite graph with very high probability.
Clearly, we cannot expect that such randomized algorithm can produce a proper coloring of the whole graph with probability one.
That is why, we need proper labeling.
\begin{definition}
Let $\cG\subset Gr_d$ be a multiset of graphs (that is each graph $G$ can be
chosen more than once). Equip each $G\in\cG$
with a Cantor labeling $\lambda_G:V(G)\to \{0,1\}^\N$. We say that the
family $\{\lambda_G\}_{G\in\cG}:V(G)\to \{0,1\}^\N$ is a proper labeling if
for any $r\geq 1$, there exists $s\geq 1$ such that
if $G\in \cG$ and $x,y\in V(G), 0<d_G(x,y)\leq r$, then the rooted-labeled balls $B_s(G,x,(\lambda_G)_s)$ and $B_s(G,y,(\lambda_G)_s)$ 
are not rooted-labeled isomorphic. 
\end{definition}
\noindent
Now, we can formally define our preferred notion of local algorithm.
\begin{definition}
Let $\cG\subset Gr_d$ be a multiset of countable graphs as above and $\{\lambda_G\}_{G\in\cG}$ be a proper labeling.
Let $m>0$ and $Q$ be a finite set. Also, let $\Theta: CU^m_d\to Q$.
Then the $Q$-labeling algorithm associated to the oracle $\Theta$ is defined
by the function $\lambda^\Theta_G:V(G)\to Q$ for each $G\in \cG$ by the rule
$$\lambda^\Theta_G(x)=\Theta(B_m(G,x,(\lambda_G)_m)).$$
\end{definition}
\noindent
Note that it is not very hard to construct proper labelings.
In \cite{Elekurs} it was proved that there exists a proper labeling
$\{\lambda_G\}_{G\in \cG}$ such that $\lambda_G:V(G)\to \{0,1\}^l$ for some $l\geq 1$. So, one can break
the symmetries using labelings by large enough finite sets.
On the other hand, it is easy to construct a proper labeling
$$\{\lambda_G:V(G)\to \{0,1\}^N\}_{G\in Gr_d}$$
\noindent
such that for any $r\geq 1$, there exists $s\geq 1$ so that for any $G\in Gr_d$ and $x,y\in V(G), 0<d_G(x,y)\leq r$ we have that
$$(\lambda_G)_s(x)\neq (\lambda_G)_s(y)\,.$$
Note that if the proper graph $E$ is a labeled realization of a Cantor labeled graph $(G,\lambda)$, then $G$ is properly labeled. On the other hand, if $(G,\lambda)$ is 
properly labeled it is not necessarily
true that a labeled realization of $(G,\lambda)$ is a proper Cantor graph, this is due to the fact the for graphs the definition of proper labeling involves
the geometric structure. \vi 
{\bf 7. Local partitioning oracles and almost finiteness.} \vi
Hassidim, Kelner, Nguyen and Onak \cite{HKNO} constructed a randomized local graph partitioning algorithm that takes finite planar graphs
(or any other hyperfinite family of graphs) of bounded vertex degrees
as an input and cuts them into small pieces by removing only a small amount of edges. 
This construction can be used
e.g. to approximate maximal independent sets in planar graphs with high probability in constant time. 
The topological analogue of the Hassidim-Kelner-Nguyen-Onak algorithm is a local algorithm that performs the task of almost finite partitionings
of finite or countable graphs. Let $\cG\subset Gr_d$ be a multiset of graphs as in the
previous section, equipped with a proper labeling
$\{\lambda_G:V(G)\to \{0,1\}^{\N}\}_{G\in \cG}$. 
\begin{definition}
Let $\epsilon>0$ be
a real number and $K\geq 1$, $m\geq 1$ be integers and $Q$ be a finite set.
Then, an $(\eps, K)$-{\bf local partitioning oracle} is given by
a function $\Theta:CU^m_d\to Q$ satisfying the following conditions.
\begin{itemize}
\item For the resulting labelings $\{\lambda_G^\Theta:V(G)\to Q\}_{G\in \cG}$,
the equality $\lambda^\Theta_G(x)=\lambda^\Theta_G(y)$ implies that
either $d_G(x,y)\leq K$ or $d_G(x,y)\geq 3K$.
Observe that this condition means that $\lambda^\Theta_G$ defines an
equivalence relation on the set $V(G)$, where $x\equiv_\Theta^\lambda y$
if $d_G(x,y)\leq K$ and $\lambda^\theta_G(x)=\lambda^\theta_G(y)$. Clearly,
the diameter of the equivalence classes are less or equal than $K$.
\item For any equivalence class $L$ and for any $G\in \cG$ we have that $i_G(L)\leq \eps$.
\end{itemize}
\end{definition}
\noindent
The next theorem is the algorithm-theoretical version of Theorem \ref{fotetel}.
\begin{theorem} \label{algo}
Let $D$ be an integer and let $\cG\subset Gr_d$ be a multiset equipped with a proper labeling
$\{\lambda_G:V(G)\to \{0,1\}^{\N}\}_{G\in \cG}$ such
that for any $G\in\cG$, $G$ is a $D$-doubling graph.  Then
for any $\eps>0$, there exist integers $m,K\geq 1$, a finite set $Q$ and
an $(\eps,K)$-local partitioning oracle $\Theta:CU^m_d\to Q$.
\end{theorem}
\noindent
In order to state our next result, we recall the notion of a {\bf local verifier}.
Let $l\geq 1$, $Q$ be a finite set and $\Omega:U^{l,Q}_d\to\{Y,N\}$ be an oracle function.
Here, $U^{l,Q}_d$ denotes the isomorphism classes of rooted $Q$-labeled $l$-balls of vertex degree bound $d$.
Now let $G\in Gr_d$ be a graph and $\phi:V(G)\to Q$ be a labeling. We say that the verifier $\Omega$
accepts the labeled graph $(G,\phi)$ if for all $x\in V(G)$, $\Omega(B_l(G,x,\phi))=Y$.
Let $\cal{H}\subset Gr_d$ be a class of finite graphs. We say that
$\Omega:U^{l,Q}_d\to\{Y,N\}$ is a local verifier for the $\eps$-approximating maximal independent
set problem in the class $\cal{H}$, if for any $G\in\cal{H}$ and for any labeling $\phi_G:V(G)\to\{0,1\}$
the verifier $\Omega$ accepts $(G,\phi_G)$ only if 
\begin{itemize}
\item $(\phi_G^{-1})(0)$ is an independent set of $G$.
\item $|(\phi_G^{-1})(0)|\geq (1-\eps) I(G)$, where $I(G)$ is the size of the largest
independent set in $G$.
\end{itemize}
\noindent
Notice that we have ``only if'' and not ``if and only if''. 
Now we can state the full topological version of the constant-time approximation result of \cite{HKNO}.
\begin{theorem} Let $D$ be an integer and let $\cG\subset Gr_d$ be a multiset of finite graphs equipped with a proper labeling
$\{\lambda_G:V(G)\to \{0,1\}^{\N}\}_{G\in \cG}$ such
that for any $G\in\cG$, $G$ is a $D$-doubling graph.  Then
for any $\eps>0$, there exist integers $m,l\geq 1$ and oracles $\Theta_\eps:CU^m_d\to \{0,1\}$,
$\Omega_\eps:U^{l,Q}_d\to\{Y,N\}$
such that
\begin{itemize}
\item $\Omega_\eps$ is a local verifier for the $\eps$-approximating maximal independent
set problem in the class of finite $D$-doubling graphs.
\item For any $G\in \cG$, $(\lambda_G^{\Theta_\eps})^{-1}(0)$ is an independent set,
\item $|(\lambda_G^{\Theta_\eps})^{-1}(0)|\geq (1-\eps) I(G)$,
\item and  $\Omega_\eps$ accepts $\lambda_G^{\Theta_\eps}$.
\end{itemize}
\end{theorem}
\vi
{\bf 8. Spectral convergence.} \vi
For a graph $G$ of bounded vertex degrees one can consider its graph Laplacian $\Delta_G:L^2(V(G))\to L^2(V(G))$ and its
spectrum $\spec(G)$ a compact subset of the non-negative real numbers. It is well-known that there exist isospectral pairs of finite graphs
that is, non-isomorphic connected finite graphs $G$ and $H$, such that $|V(G)|=|V(H)|$ and $\spec(G)=\spec(H)$ see e.g. \cite{GM} for
lots of examples. For infinite graphs, much less is known besides that isospectral pairs do exist
\cite{kinai}. However, we will show that
if two connected infinite graphs $G$ and $H$ are equivalent, then $\spec(G)=\spec(H)$. This shows that one cannot ``hear''  the shape
of an infinite graph, not even in an approximate sense.
\begin{prop}
There exist uncountably many pairwise non-quasi isomorphic graphs that share the same spectrum.
\end{prop}
\noindent
If $G$ is a finite graph then
there is also a natural probability measure $\mu_G$ on the spectrum of $\Delta_G$, when $\mu_G(\lambda)$ is the multiplicity of $\lambda$ in the spectrum
divided by $|V(G)|$. 
Let $\{G_n\}^\infty_{n=1}$ be a sequence of finite graphs that are convergent in the sense of Benjamini and Schramm. It is easy to see that the
sequence of Borel measures $\{\mu_{G_n}\}^\infty_{n=1}$ are weakly convergent. One can ask, what happens in the case of topological convergence. Using our result
on fractional almost finiteness we can prove the following theorem.
\begin{theorem} \label{specconv} Let $D>0$ and $\{G_n\}^\infty_{n=1}$ be a sequence of $D$-doubling countable 
graphs converging to a graph $G$ in the topological sense. 
Then $\spec(G_n)\to \spec(G)$ in the Hausdorff topology. \end{theorem}
\noindent
Note that the theorem above holds for the spectra of arbitrary local kernel operators (generalized Laplacians)
as well. Results on Hausdorff convergence of the spectra of  combinatorial operators on infinite
graphs can also be found in \cite{BB}. 
Now let us consider the Cayley graph of a finitely generated non-amenable group $\Gamma$ with respect to a symmetric generating system $\Sigma$. 
Let $\Gamma\supset N_1 \supset N_2,\dots$ be
a nested sequence of finite index normal subgroups such that
$\cap_{k=1}^\infty{N_k}=\{1\}$. Then the sequence of finite Cayley graphs $\{\mbox{Cay}(\Gamma/N_k,\Sigma)\}^\infty_{k=1}$ topologically converges to the Cayley graph
$\mbox{Cay}(\Gamma,\Sigma)$. Nevertheless, the spectra of the finite Cayley graphs contain the zero and the spectrum of $\mbox{Cay}(\Gamma,\Sigma)$ is away from
the zero. Hence, the spectra of the finite Cayley graphs do not converge to the spectrum of the Cayley graph of $\Gamma$ in the Hausdorff topology.
However, we have the following theorem.
\begin{theorem} \label{luck}
Let $\Gamma$ be a finitely generated amenable group with symmetric generating system $\Sigma$ and let $\Gamma\supset N_1 \supset N_2,\dots$ be
a nested sequence of finite index normal subgroups such that $\cap_{k=1}^\infty{N_k}=\{1\}$. Then the spectra of the Cayley graphs $\{\mbox{Cay}(\Gamma/N_k,\Sigma)\}^\infty_{k=1}$ converge to the spectrum of $\mbox{Cay}(\Gamma,\Sigma)$ in the Hausdorff topology. \end{theorem}
\vi
{\bf 9. $C^*$-algebras.}
\vi
For free continuous Cantor actions $\alpha:\Gamma\actson \cC$ one can consider its reduced $C^*$-algebras. Recently, Kerr \cite{Kerr} was able to prove
that almost finiteness of a free, minimal action implies the so-called $\cal{Z}$-stability of the associated reduced $C^*$-algebras.
By the result above,  Theorem \ref{fotetel} implies that the free, minimal actions of finitely generated nilpotent-by-finite groups has $\cal{Z}$-stable
reduced $C^*$-algebra. 
In \cite{Suzuki}, Suzuki proved that the groupoid $C^*$-algebras of minimal, almost finite, \'etale groupoids have stable rank one, that is, the invertible
elements form a dense subset. In \cite{Elekurs}, we associated certain $C^*$-algebras to uniform recurrent subgroups and even arbitrary proper Cantor actions
, using
local kernel operators. Aaron Tikuisis pointed out to me that
in the case of proper actions these are exactly the groupoid algebras and the corresponding groupoid is \'etale. Hence, combining Suzuki's result with
the construction of various proper minimal actions in \cite{Elekurs} one can find many non-free actions that result in simple $C^*$-algebras of stable rank one.

\end{document}